\documentclass[leqno]{article}
\usepackage{amsmath,amsfonts}

\numberwithin{equation}{section}

\newtheorem{prop}{Proposition}[section]
\newtheorem{defn}[prop]{Definition}

\newtheorem{thm}[prop]{Theorem}

\newtheorem{lem}[prop]{Lemma}

\newcommand\cI{{\mathcal{I}}}

\newcommand\nr{N_{\mathbb{R}}}
\newcommand\R{\mathbb{R}}
\newcommand\Z{\mathbb{Z}}
\newcommand\C{\mathbb{C}}

\newcommand\p{\mathbb{P}}

\newcommand\cO{{\mathcal{O}}}

\newcommand\prf{\noindent {\bf Proof.\ \ }}
\newcommand\qed{\hfill$\clubsuit$}

\begin{document}

\title{Equations of Parametric Curves and Surfaces via Syzygies}

\author{David A.\ Cox\\
Department of Mathematics and Computer Science\\
Amherst College\\
Amherst, MA 01002-5000\\
{\tt dac@cs.amherst.edu}}



\maketitle

\section{Implicitization}

Suppose that $\phi : \p^2 \to \p^3$ is a map (possibly with
basepoints) whose image is a surface $S \subset \p^3$.  In computer
graphics, such a map is given by homogeneous polynomials $(a,b,c,d)$
with real coefficients, and knowing these polynomials allows one to
draw the real points of $S$ on a computer screen.  Most geometric
models in CAD (Computer Aided Design) use parametric surfaces.  This
includes car bodies, airplanes, and animated figures in movies such as
\emph{Toy Story} or \emph{Dinosaur}. 

The \emph{implicitization problem}, as explained in \cite{clo1}, is to
compute the implicit equation $F = 0$ of $S$ using the parametrization
$\phi$.  This leads to the following natural question:\ if our goal is
to draw $S$ on a computer screen, why would we be interested in the
implicit equation of $S$?

One answer is that we can use implicitization to help find curve
intersections.  Namely, when drawing two parametrized surfaces in
3-dimensional space, one sometimes wants to highlight the curve
where they intersect.  Suppose that we have parametrizations $\phi_i :
\p^2 \to \p^3$, $i=1,2$, and for simplicity assume that there are no
basepoints.  Then we want to find $S_1\cap S_2$, where $S_i$ is the
image of $\phi_i$.

To solve this problem, let $F_1 = 0$ be the implicit equation of
$S_1$, and let $C \subset \p^2$ be defined by $F_1 \circ \phi_2 = 0$.
This curve may be singular of high genus, but since we have an
explicit equation for $C$, there are known algorithms for drawing it
(see, for example, \cite[Section 6.5]{hoffmann}).  Then $\phi_2(C)$ is
the desired intersection $S_1\cap S_2$.

Another use of implicitization occurs when one creates new geometric
models by applying Boolean operations to existing models.  For
example, if two balls interesect in $\R^3$, then removing one from the
other involves knowing the curve where their boundaries intersect.

The next question concerns how to find the implicit equation, assuming
we are given the parametrization $\phi$.  In practice, three methods
are used:
\begin{itemize}
\item Gr\"obner bases.
\item Resultants.
\item Syzygies.
\end{itemize}
This paper will concentrate on the third of these methods.

\section{Syzygies and Equations of Curves}

This section will report on joint work \cite{csc} with Tom Sederberg
(Brigham Young University) and Falai Chen (University of Science and
Technology of China). 

We will work over $\C$.  The idea is that we want to implicitize the
map
\[
\phi : \p^1 \longrightarrow \p^2,
\]
given by 
\begin{equation}
\label{curve.1}
\phi(s,t) = \big(a(s,t),b(s,t),c(s,t)\big),
\end{equation}
where $a,b,c \in R = \C[s,t]$ are homogeneous polynomials of degree
$n$ (here, $s,t$ are homogeneous coordinates on $\p^1$).  We will also
assume that $\gcd(a,b,c) = 1$.  This ensures that $\phi$ has no
basepoints (a point $(s_0,t_0) \in \p^1$ is a \emph{basepoint} of
$\phi$ if $a,b,c$ all vanish at $(s_0,t_0)$).

In \cite{ssqk,sc,sgd}, Sederberg and his co-workers introduced the
idea of a \emph{moving line} in $\p^2$.  If we let $x,y,z$ be
homogeneous coordinates for $\p^2$, then a moving line is an equation
of the form
\begin{equation}
\label{moving.1}
A(s,t) x + B(s,t) y + C(s,t) z = 0,
\end{equation}
where $A,B,C \in R$ are homogeneous of the same degree.  We can regard 
\eqref{moving.1} as a family of lines parameterized by $(s,t) \in
\p^1$.  

One can easily imagine how the point of intersection of two moving
lines traces out a curve in $\p^2$.  This leads to the question of
whether the map $\phi$ from \eqref{curve.1} arises this way.  If we
dehomogenize by setting $t = 1$ (in $\p^1$) and $z = 1$ (in $\p^2$),
then we get an easy answer, for in the case, \eqref{curve.1} gives the
curve in $\C^2$ parametrized by
\[
\begin{aligned}
x &= \frac{a(s)}{c(s)}\\[2pt]
y &= \frac{b(s)}{c(s)}.
\end{aligned}
\]
where $a(s)$ is short for $a(s,1)$, and similarly for $b(s)$ and
$c(s)$.  This can be thought of as the point of intersection of the
moving vertical line $c(s)x - a(s) = 0$ and the moving horizontal line
$c(s)y - b(s) = 0$.  Note that these are moving lines of degree $n$.
As we will see below, we get significantly lower degrees by allowing
more general moving lines.

To formalize the above discussion, we make the following definition.

\begin{defn} 
\label{follow.1}
The moving line \eqref{moving.1} {\bf follows} the parametrization
\eqref{curve.1} if
\[
A(s,t) a(s,t) + B(s,t) b(s,t) + C(s,t) c(s,t) = 0
\]
for all $(s,t) \in \p^1$.
\end{defn}

Geometrically, this means that for all $(s,t)$, the point on the
parametrized curve lies on the corresponding line.  More important is
the algebraic interpretation of Definition~\ref{follow.1}, which says
that $A,B,C$ is a \emph{syzygy} on $a,b,c$.  We write this as
\[
(A,B,C) \in \mathrm{Syz}(a,b,c),
\]
where $\mathrm{Syz}(a,b,c) \subset R^3$ is the syzygy module of
$(a,b,c)$. 

Since $\mathrm{Syz}(a,b,c)$ is a graded module, we can speak of its
graded piece in dimension $d$, denoted $\mathrm{Syz}(a,b,c)_d$.  We
will now describe how $\mathrm{Syz}(a,b,c)_{n-1}$ determines the
implicit equation of the image of \eqref{curve.1}.

To see how this works, consider the map
\begin{equation}
\label{ml.1}
R_{n-1}^3 \xrightarrow{(a,b,c)} R_{2n-1},
\end{equation}
where subscripts indicated graded pieces (remember that $a,b,c$ have
degree $n$).  The kernel of this map is $\mathrm{Syz}(a,b,c)_{n-1}$.
Note also that $\dim R_{n-1}^3 = 3n$ and $\dim R_{2n-1} = 2n$, so that
\begin{equation}
\label{maximal.1}
\dim \mathrm{Syz}(a,b,c)_{n-1} = n \iff \text{\eqref{ml.1} has maximal
rank.} 
\end{equation} 

Later, we will see that \eqref{ml.1} always has maximal rank by
regularity.  We will assume this for now.  By \eqref{maximal.1}, it
follows that we can find $n$ linearly independent moving lines which
follow $\phi$.  Write these moving lines as follows:
\begin{equation}
\label{lij.1}
A_i x + B_i y + C_i z = \sum_{j=0}^{n-1} L_{i,j}(x,y,z)\,s^j
t^{n-1-j},\quad i = 0,\dots,n-1.
\end{equation}
Note that $L_{i,j}(x,y,z)$ is a linear form with coefficients in $\C$.
Then one of the main results of \cite{sc} is the following.

\begin{thm}
\label{impl.1}
Let $C$ be the image of \eqref{curve.1}, and let $d$ be the generic
degree of the induced map $\p^1 \to C$.  Then
\[
\det (L_{i,j}) = \lambda\,F^d,
\]
where $\lambda \in \C\setminus\{0\}$ and $F = 0$ is the
$($irreducible\/$)$ implicit equation of $C \subset \p^2$.
\end{thm}

Note how the numbers work in this theorem: since $a,b,c$ have degree
$n$, the curve $C$ traced out by $\phi$ has degree $n/d$, where $d$ is
the generic degree.  Thus $F^d$ has degree $n$.  On the other hand,
the determinant in Theorem~\ref{impl.1} also has degree $n$ since
$(L_{i,j})$ is an $n\times n$ matrix of linear forms.

So far, we have used only one graded piece of the syzygy module,
namely $\mathrm{Syz}(a,b,c)_{n-1}$.  We next turn our attention to the
entire syzygy module.  As we will see below, the structure of this
module will give deeper insight into the determinant $\det (L_{i,j})$
used in Theorem~\ref{impl.1}.

The key tool for understanding $\mathrm{Syz}(a,b,c)$ is the Hilbert
Syzygy Theorem, which implies that syzgy modules of homogeneous
polynomials in $R = \C[s,t]$ are always free.  More precisely, if we
set $I = \langle a,b,c\rangle \subset R = \C[s,t]$, then the Syzygy
Theorem and an easy argument using the Hilbert Polynomial imply that
$R/I$ has a free resolution
\[
0 \to R(-n-\mu_1) \oplus R(-n-\mu_2) \to R(-n)^3 \to R \to R/I \to 0,
\]
where $\mu_1+\mu_2 = n$ and the map $R^3 \to R$ is given by $a,b,c$.
(See \cite{clo2,csc} for the details.)  We can assume $\mu_1 \le
\mu_2$, and following \cite{csc}, we let $\mu = \mu_1$, so that $\mu
\le n-\mu = \mu_2$.  In down to earth terms, the above resolution
means that $\mathrm{Syz}(a,b,c)$ is a free $R$-module of rank two with
generators, say $p$ and $q$, of respective degrees $\mu$ and $n-\mu$.
We call $p,q$ a $\mu$-\emph{basis} of the syzygy module.

The existence of a $\mu$-basis has some nice consequences.  First,
note that syzygies of degree $n-1$ can be uniquely written in the form
\begin{equation}
\label{nm1syz.1}
h_1\,p + h_2\,q,
\end{equation}
where
\[
\deg h_1 = n-\mu-1\quad\text{and}\quad \deg h_2 =
\mu-1.
\]
Computing dimensions, we conclude that 
\[
\dim \mathrm{Syz}(a,b,c)_{n-1} = n, 
\]
so that by \eqref{maximal.1}, we see that \eqref{ml.1} always
has rank $n$, as claimed earlier.

Another consequence \eqref{nm1syz.1} is that a basis of
$\mathrm{Syz}(a,b,c)_{n-1}$ is given by
\[
s^it^{n-\mu-1-i} p,\ 0 \le i \le n-\mu-1,\quad s^it^{\mu-1-i} q,\ 0
\le i \le \mu-1.
\]
If we use this basis to form the linear forms $L_{i,j}$ as in
\eqref{lij.1}, then we easily obtain
\[
\det (L_{i,j}) = \mathrm{Res}(p,q).
\]
Thus the syzygies of degree $n-1$ compute the resultant of the
$\mu$-basis.  If we combine this with Theorem~\ref{impl.1}, then we
conclude that
\[
\mathrm{Res}(p,q) = \lambda\,F^d,\ \lambda \in \C\setminus\{0\}.
\]
This shows that the $\mu$-basis computes the implicit equation of the
curve.  A careful proof can be found in \cite{csc}.  

We can also make some comments from a computational point of view.  In
order to use Theorem~\ref{impl.1}, one needs a basis of
$\mathrm{Syz}(a,b,c)_{n-1}$.  Since the matrix of \eqref{ml.1} is
$2n\times3n$, the complexity of computing a basis is $O(n^3)$ by
standard methods in numerical linear algebra.  However, \cite{sz}
gives an $O(n^2)$ method, based on variant of the Buchberger
algorithm, for finding a $\mu$-basis of $\mathrm{Syz}(a,b,c)$.

Here are some other results from \cite{csc} concerning $\mu$-bases:
\begin{itemize}
\item The regularity of the ideal $I = \langle a,b,c\rangle \subset R
= \C[s,t]$ is $n-\mu-1$.  Thus knowing $\mu$ is equivalent to knowing
the regularity.  Note also that $n-\mu-1 \le n-1$.  This explains why
syzygies of degree $n-1$ work so well---regularity always holds for
this degree, no matter what $\mu$ is.
\item We defined $\mu$ so that $\mu \le n-\mu$, which implies $0
\le \mu \le \lfloor n/2\rfloor$.  When one considers all triples
$a,b,c$ in $R$ of degree $n$ with $\gcd(a,b,c) = 1$, one can show that
$\mu = \lfloor n/2\rfloor$ is generic.  
\item The usual Sylvester form of the resultant expresses
$\mathrm{Res}(p,q)$ as the $n\times n$ determinant described above.
One can also express $\mathrm{Res}(p,q)$ as a $(n-\mu)\times (n-\mu)$
determinant where $n-2\mu$ rows are linear in the entries of $p$ and
$\mu$ rows are quadratic and are built from the B\'ezoutian
of $p$ and $q$.  
\end{itemize}

Finally, I should comment that ``moving lines'' represent an
independent discovery of the concept of syzygy by the computer science
community.  Furthermore, in their definition of ``$\mu$-basis'',
Sederberg and Chen essentially conjectured a special case of the
Hilbert Syzygy Theorem.  

It is interesting to note that this special case, which asserts that
$\mathrm{Syz}(a,b,c)$ is a graded free $R$-module, was actually proved
by Franz Meyer in 1887 in \cite{meyer}.  Meyer also conjectured that a
similar result should hold for $\mathrm{Syz}(a_1,\dots,a_m)$, where
the $a_i \in R = \C[s,t]$ are homogeneous, but he was unable to prove
this in general.  Hilbert, in his great paper \cite{hilbert} of 1890,
proves the general form of the Hilbert Syzygy Theorem and explains how
to use the Hilbert polynomial.  His very first application is to
Meyer's conjecture from 1887.

\section{Syzygies and Tensor Product Surfaces} 

We now turn our attention to surface parametrizations.  This section
will report on the paper \cite{cgz} written with Ron Goldman and Ming
Zhang, both of Rice University.  Again, most proofs will be omitted.

In this section, we will consider a \emph{tensor product
parametrization}, which is a map
\[
\phi : \p^1\times\p^1 \longrightarrow \p^3
\]
defined by
\begin{equation}
\label{surface.3}
\phi(s,u;t,v) = \big(a(s,u;t,v),b(s,u;t,v),c(s,u;t,v),d(s,u;t,v)\big),
\end{equation}
where $a,b,c,d \in R = \C[s,u;t,v]$ are bihomogeneous polynomials of
bidegree $(m,n)$.  Here, we think of $s,u$ as homogeneous coordinates
on the first factor of $\p^1$ and $t,v$ as homogeneous coordinates on
the second.  We will also assume that $a,b,c,d$ have no common
factors, which implies that $\phi$ has at most finitely many
basepoints (a point $(s_0,u_0;t_0,v_0) \in \p^1\times\p^1$ is a
\emph{basepoint} if $a,b,c,d$ all vanish at $(s_0,u_0;t_0,v_0)$). 

If $\phi$ has no basepoints (which we will assume throughout this
section) and is generically one-to-one, then it is well-known that the
image of $\phi$ is a surface $S \subset \p^3$ of degree $2mn$.  The
goal of this section is to find the defining equation of $S$ using
syzygies. 

The analog of a moving line in $\p^2$ is clearly a \emph{moving plane}
in $\p^3$.  This is an equation of the form
\[
A(s,u;t,v) x + B(s,u;t,v)y + C(s,u;t,v)z + D(s,u;t,v)w = 0,
\]
where $x,y,z,w$ are homogeneous coordinates on $\p^3$ and $A,B,C,D \in
R$ are bihomogeneous of the same bidgree.  Moving planes were first
considered in \cite{sc}.  Then we say that the above
moving plane \emph{follows} the parametrization \eqref{surface.3} if
\begin{align*}
A(s,u;t,v)a(s,u;t,v) +\ &B(s,u;t,v)b(s,u;t,v) + \\
&C(s,u;t,v)c(s,u;t,v) +
D(s,u;t,v)d(s,u;t,v) = 0
\end{align*}
for all $(s,u;t,v) \in \p^1\times\p^1$.  Thus the moving
plane follows the parametrization if and only if
\[
(A,B,C,D) \in \mathrm{Syz}(a,b,c,d).
\]

As we will soon see, moving planes are not sufficient---to get the
implicit equation of the image of $\phi$, we will need to use moving
surfaces of higher degree.  This idea appears in \cite{sc}, and even
for the curve case, one can use \emph{moving conics} to get some
interesting results concerning implicitization (see \cite{sc,zcg}).

For us, the crucial ingredient will be \emph{moving quadrics}, which
are equations of the form
\[
A(s,u;t,v) x^2 + B(s,u;t,v) xy + \dots + I(s,u;t,v) zw + J(s,u;t,v)
w^2 = 0,
\]
where $A,B,\dots,I,J \in R$ are bihomogeneous of the same bidgree.
It should be clear what it means for a moving quadric to
\emph{follow} the parametrization \eqref{surface.3}, and one easily
sees that this is equivalent to 
\[
(A,B,\dots,I,J) \in \mathrm{Syz}(a^2,ab,\dots,cd,d^2) \subset R^{10}.
\]

For curves, moving lines of degree $n-1$ played a crucial role.  For a
tensor product surface, it thus makes sense to consider moving planes
and quadrics of bidegree $(m-1,n-1)$ which follow the parametrization.
The moving planes of this degree are the kernel of the map
\begin{equation}
\label{mp.3}
MP : R_{m-1,n-1}^4 \xrightarrow{(a,b,c,d)} R_{2m-1,2n-1}
\end{equation}
(remember that $a,b,c,d$ have bidegree $(m,n)$).  Both of these vector
spaces have dimension $4mn$, so that generically, we expect $MP$ to be
an isomorphism.  In other words, there should usually be \emph{no}
moving planes of bidegree $(m-1,n-1)$.  

Thus we turn our attention to moving quadrics which follow the
parametrization.  In bidegree $(m-1,n-1)$, these are given by the
kernel of the map
\begin{equation}
\label{mq.3}
MQ : R_{m-1,n-1}^{10} \xrightarrow{(a^2,ab,\dots,cd,d^2)}
R_{3m-1,3n-1} 
\end{equation}
In this case, one easily sees that 
\[
\dim \mathrm{Syz}(a,b,c,d)_{m-1,n-1} = mn \iff \text{\eqref{mq.3} has
maximal rank.} 
\]

For now, we will assume that we have precisely $mn$ linearly
independent moving quadrics of bidegree $(m-1,n-1)$ which follow the
parametrization.  Label these as $Q_i$ for $1 \le i \le mn$.  The idea
is to construct a square matrix by writing out the $Q_i$ as we did in
\eqref{lij.1}.  Here, we first dehomogenize by setting $u = v = 1$
to simplify the resulting formulas.  Thus $Q_i$ can be written
\begin{equation}
\label{expand.3}
\begin{aligned}
Q_i &= A_i x^2 + \dots + J_i w^2\\[-1pt]
	&= \Big(\sum_{j=0}^{m-1} \sum_{k=0}^{n-1} A_{i,jk} s^jt^k\Big)
	x^2 	+ \dots + \Big(\sum_{j=0}^{m-1} \sum_{k=0}^{n-1}
	J_{i,jk} s^jt^k\Big) w^2\\[-1pt]
	&= \sum_{j=0}^{m-1} \sum_{k=0}^{n-1} \Big(A_{i,jk} x^2 + \dots
	+ J_{i,jk} w^2\Big)s^jt^k\\[-1pt]
	&= \sum_{j=0}^{m-1} \sum_{k=0}^{n-1} Q_{i,jk}(x,y,z,w)s^jt^k,
\end{aligned}
\end{equation}
where $Q_{i,jk}$ is a quadric in $x,y,z,w$ with coefficients in $\C$.
Furthermore, note that $i$ ranges over the $mn$ numbers 1 to $mn$
while $(j,k)$ ranges over the $mn$ pairs $(0,0)$ to $(m-1,n-1)$.  It
follows that we can arrange the $Q_{i,jk}$ into a square matrix of
size $mn\times mn$, where each entry is a quadric in $x,y,z,w$.  We
write this as
\[
M = (Q_{i,jk}).
\]
Notice that $\det M$ has degree $2mn$ in $x,y,z,w$.

We can now state one of the main results of \cite{cgz}.

\begin{thm}
\label{tensor.3}
Suppose that $\phi : \p^1\times\p^1 \to \p^3$ has no basepoints and is
generically one-to-one.  If $MP$ from \eqref{mp.3} has maximal rank,
then so does $MQ$ from \eqref{mq.3} and furthermore, the image of
$\phi$ is defined by the equation $\det M = 0$.
\end{thm}

\prf I will sketch some parts of the proof.  One begins by
changing coordinates on $\p^3$ so that $a,b,c$ have no basepoints.
Then consider the map
\[
MQ' : R_{m-1,n-1}^{9} \xrightarrow{(a^2,ab,\dots,cd)}
R_{3m-1,3n-1}
\]
obtained from \eqref{mq.3} by omitting $d^2$.  Here, $\dim
R_{m-1,n-1}^{9} = \dim R_{3m-1,3n-1} = 9mn$.  Thus $\det MQ' \ne 0$
implies that $MQ$ has maximal rank.  This in turn will give the $mn$
linearly independent moving quadrics of bidegree $(m-1,n-1)$ needed to
construct the matrix $M$.

I will discuss two proofs that $\det MQ' \ne 0$.  For the first,
suppose that $\det MQ' = 0$.  Then there is a nontrivial syzygy
\[
Aa^2 + Bab + \dots + I cd = 0,
\]
where $A,B,\dots,I$ are bihomogeneous of bidegree $(m-1,n-1)$.  Since
every term contains $a$, $b$ or $c$ (we got rid of $d^2$), we obtain
\begin{equation}
\label{abcsyz.3}
(Aa+Bb+Cc+Dd)a + (Eb+Fc+Gd)b + (Hc+Id)c = 0.
\end{equation}
This is a syzygy on $a,b,c$ of bidegree $(2m-1,2n-1)$.  I remember
when Ron Goldman showed me this equation and asked me if it implied
that
\begin{equation}
\label{abckoszul.3}
Hc + Id = -h_1a - h_3 b
\end{equation}
for bihomogeneous polynomials $h_1,h_3$ of bidegree
$(m-1,n-1)$.  If \eqref{abckoszul.3} is true, then we get a nontrivial
syzygy on $a,b,c,d$, which contradicts our assumption that $MP$ has
maximal rank.

Hence I needed to show that \eqref{abcsyz.3} implies
\eqref{abckoszul.3}.  The idea, of course, is that there is a Koszul
complex lurking in the background.  In general, if
\[
\mathcal{A}a + \mathcal{B}b + \mathcal{C}c = 0,
\]
then we say that $(\mathcal{A},\mathcal{B},\mathcal{C})$ is a
\emph{Koszul syzygy} if there are $h_1,h_2,h_3$ such that
\begin{equation}
\label{koszulsyz.3}
\begin{aligned}
	\mathcal{A} &= \ h_1c + h_2b\\
	\mathcal{B} &= -h_2a + h_3c\\
	\mathcal{C} &= -h_1a - h_3b.
\end{aligned}
\end{equation}

If we were working in $\p^2$, then $a,b,c$ having no basepoints would
imply that they were a regular sequence, and it would follow
immediately that their Koszul complex is exact.  In particular, every
syzygy on $a,b,c$ would be Koszul, so that \eqref{abcsyz.3}
$\Rightarrow$ \eqref{abckoszul.3} is automatic in $\p^2$.

But we are in $\p^1\times\p^1$, where $a,b,c$ are bihomogeneous.  In
this bigraded situation, $a,b,c$ almost never form a regular sequence,
and their Koszul complex need not be exact (it is easy to give
counterexamples).  Instead, I had to use the vanishing of a certain
sheaf cohomology group to show that every syzygy of bidegree
$(2m-1,2n-1)$ is Koszul (see \cite{cgz} for details).  Hence
\eqref{abckoszul.3} is true, and as explained above, this completes
the first proof that $\det MQ' \ne 0$.

The second proof that $\det MQ' \ne 0$ is much quicker.  In
\cite{zcg}, it was conjectured that 
\begin{equation}
\label{identity.3}
\det MQ' = \mathrm{Res}(a,b,c) \big(\det MP\big)^3.
\end{equation}
Since $a,b,c$ have no basepoints, their resultant is nonvanishing,
and $\det MP \ne 0$ by assumption.  Then $\det MQ' \ne 0$ would follow
immediately from the above identity.  In a recent paper
\cite{dandrea}, Carlos D'Andrea not only proved \eqref{identity.3} but
also generalized it moving surfaces of degree $> 2$ as well.

Once we know that $\det MQ' \ne 0$, we can construct the desired
matrix $M$.  It is easy to see that $\det M$ vanishes on the image of
$\phi$ since the moving quadrics used in $M$ all follow the
parametrization.  Furthermore, standard techniques from resultant
theory show that $\det M$ is not identically zero (one shows that the
coefficient of $w^{2mn}$ is nonzero).  It follows that $\det M$ is a
nonzero polynomial of degree $2mn$ which vanishes on the image of
$\phi$.  But since $\phi$ is generically one-to-one, its image has
degree $2mn$.  Hence $\det M = 0$ must be the irreducible equation of
the image.\qed

\medskip

Theorem \ref{tensor.3} assumes that $\phi$ has no basepoints, that $MP$
has maximal rank, and that $\phi$ is generically one-to-one.
Recently, D'Andrea has been able to weaken some of these hypotheses:
\begin{itemize}
\item If $MP$ has maximal rank and $d$ is the generic degree of
$\phi$, then $\det M = \lambda\,F^d$, where $\lambda \in
\C\setminus\{0\}$ and $F = 0$ is the irreducible equation of the
image.  This is proved in \cite{dandrea}.
\item If $\phi$ is generically one-to-one, the $MP$ has maximal rank.
This is unpublished.  
\end{itemize}
\noindent It follows that when $\phi$ has no basepoints, we can modify
Theorem \ref{tensor.3} to assert that $\det M = \lambda\,F^d$ if
\emph{either} $MP$ has maximal rank \emph{or} $\phi$ is generically
one-to-one.

\section{Syzygies and Triangular Surfaces} 

Here, we will indicate how the results of the previous section can be
modified in the case of a \emph{triangular parametrization}
\[
\phi : \p^2 \longrightarrow \p^3,
\]
which is given by homogeneous polynomials $a,b,c,d \in R = \C[s,t,u]$
of degree $n$.  As above, we assume that $a,b,c,d$ have no basepoints.
Then the analogs of $MP$ and $MQ$ are
\begin{equation}
\label{mpmq.4}
\begin{aligned}
MP &: R_{n-1}^4 \xrightarrow{(a,b,c,d)} R_{2n-1}\\[2pt]
MQ &: R_{n-1}^{10} \xrightarrow{(a^2,ab,\dots,cd,d^2)} R_{3n-1}.
\end{aligned}
\end{equation}
Assuming that $MP$ and $MQ$ have maximal rank, one easily computes
that 
\[
\begin{aligned}
\dim \mathrm{Syz}(a,b,c,d)_{n-1} &= n\\[2pt]
\dim \mathrm{Syz}(a^2,ab,\dots,cd,d^2)_{n-1} &= (n^2+7n)/2.
\end{aligned}
\]
However, each moving plane of degree $n-1$ which follows $\phi$ can be
multiplied by $x$, $y$, $z$, or $w$ to get a moving quadric.  This
gives a subspace of $\mathrm{Syz}(a^2,\dots,d^2)_{n-1}$ of dimension
$4n$, and if we pick a complementary subspace, then we obtain
\[
\begin{aligned}
n &\quad \text{linearly independent moving planes which follow}\
\phi\\[4pt]
(n^2-n)/2 &\quad \text{linearly independent moving quadrics not coming}\\
	&\quad \text{from moving planes which follow}\ \phi.
\end{aligned}
\]
Note also that there are $(n^2+n)/2$ monomials in $s,t,u$ of degree
$n-1$.  It follows that if we expand the above moving planes and
quadrics as in \eqref{expand.3}, then we get a matrix $M$ of size
$(n^2+n)/2 \times (n^2+n)/2$, where the first $n$ rows (coming from
the moving planes) are linear in $x,y,z,w$ and the remaining
$(n^2-n)/2$ rows (coming from moving quadrics) are quadratic in
$x,y,z,w$.  It follows that
\[
\deg(\det M) = 1\cdot n + 2\cdot (n^2-n)/2 = n^2.
\]
Then the following theorem is proved in \cite{cgz}.

\begin{thm}
\label{triangular.4}
Suppose that $\phi : \p^2 \to \p^3$ has no basepoints and is
generically one-to-one.  If $MP$ from \eqref{mpmq.4} has maximal rank,
then so does $MQ$ from \eqref{mpmq.4} and furthermore, the image of
$\phi$ is defined by the equation $\det M = 0$.
\end{thm}

As in the tensor product case, the exactness of a certain Koszul
complex plays a key role in the proof.  We also note that the
improvements that D'Andrea made to Theorem~\ref{tensor.3} apply to
Theorem~\ref{triangular.4} as well.

\section{Syzygies and Basepoints} 

This section discusses new results which (we hope!)\ will shed light
on how syzygies can be used to compute implicit equations of
parametrized surfaces in the presence of basepoints.  For simplicity,
we will concentrate on the triangular case, where
\[
\phi : \p^2 \,-\!\to \p^3
\]
is the rational map given by homogeneous polynomials $a,b,c,d \in R =
\C[s,t,u]$ of degree $n$.  Note that $\phi$ is a morphism outside the
set of basepoints.

\subsection{Strong \boldmath{$\mu$}-Bases for Surfaces}  We first ask
if it ever happens that the syzygy module $\mathrm{Syz}(a,b,c,d)$ is
free.  While this always happens in the curve case, it is quite rare
for surfaces.  In the discussion which follows, we will make frequent
use of standard results in commutative algebra.  A good reference is
\cite{eisenbud}, especially Chapters 18--20.

We first consider the case when $\phi$ has no basepoints.  

\begin{prop}
\label{neverfree.5}
Let $a,b,c,d \in R = \C[s,t,u]$ be homogeneous polynomials of degree
$n$, and assume that $a,b,c,d$ have no common zeros on $\p^2$.  Then
$\mathrm{Syz}(a,b,c,d)$ is not a free $R$-module.
\end{prop}

\prf
Let $I = \langle a,b,c,d\rangle \subset R$, and let $\mathfrak{m}$
denote the maximal ideal of $R/I$.  This ring has Krull dimension 0
since there are no basepoints, and thus $\mathfrak{m}$ has codimension
0.  The usual inequality between depth and codimension implies that
$\mathfrak{m}$ has depth 0 as well.  Then the Auslander-Buchsbaum
Theorem easily implies that the projective dimension of $R/I$ is 3.

However, if $\mathrm{Syz}(a,b,c,d)$ were free, then we would get the
free resolution
\[
0 \to \mathrm{Syz}(a,b,c,d) \to R(-n)^4 \to R \to R/I \to 0,
\]
which would imply that the projective dimension is $\le 2$.\qed  
\medskip

We next consider the case when there are basepoints.   

\begin{prop}
\label{strong.5}
Let $a,b,c,d \in R = \C[s,t,u]$ be homogeneous polynomials of degree
$n$, and assume that $\gcd(a,b,c,d) = 1$.  Set $I = \langle
a,b,c,d\rangle \subset R$.  If $a,b,c,d$ have at least one common zero
in $\p^2$, then the following are equivalent:
\begin{enumerate}
\item $\mathrm{Syz}(a,b,c,d)$ is a free graded $R$-module.
\item $R/I$ has projective dimension 2.
\item $R/I$ is Cohen-Macaulay.
\item $\langle s,t,u\rangle \notin \mathrm{Ass}(R/I)$.
\item $I$ is saturated.
\end{enumerate}
\end{prop}

\prf
The equivalence of (1) and (2) follows from the proof of
Proposition~\ref{neverfree.5}, and the equivalence of (2) and (3)
follows easily from the Auslander-Buchsbaum Theorem and the definition
of Cohen-Macaulay.  The equivalence of (2) and (4) follows from
Corollary~19.10 of \cite{eisenbud} since $I$ is an ideal of
codimension 2.  (This corollary is for local rings, but as is typical,
the proof also applies to the graded situation considered here.)
Finally, $\langle s,t,u\rangle \in \mathrm{Ass}(R/I)$ if and only if
$\langle s,t,u\rangle = \mathrm{Ann}(a)$ for some $a \in R/I$.  The
latter is equivalent to $I\! : \!\langle s,t,u\rangle \ne I$, and the
equivalence of (4) and (5) follows.\qed
\medskip

The original version of this proposition gave only the equivalence of
(1), (2) and (3).  I am grateful to Hal Schenk for pointing out the
relevance of (4) and (5).

Now suppose that $\mathrm{Syz}(a,b,c,d)$ is a free graded $R$-module.
Then an easy Hilbert polynomial argument shows that we have an exact
sequence of the form
\[
0 \to R(-n-\mu_1)\oplus R(-n-\mu_2)\oplus R(-n-\mu_3) \xrightarrow{\
A\ } R(-n)^4 \xrightarrow{\ B\ } R \to R/I \to 0,
\]
where $\mu_1+\mu_2+\mu_3 = n$.  Here $B$ is the map given by
$(a,b,c,d)$, and the columns of $A$ give three syzygies $p_1,p_2,p_3$
of respective degrees $\mu_1,\mu_2,\mu_3$ which are free generators of
$\mathrm{Syz}(a,b,c,d)$.  Furthermore, as is well-known, the
Hilbert-Burch Theorem implies that $a,b,c,d$ are (up to sign) the
maximal minors of the matrix $A$.

In this situation, we say that $p_1,p_2,p_3$ form a \emph{strong
{\boldmath$\mu$}-basis}, where $\mbox{\boldmath$\mu$} =
(\mu_1,\mu_2,\mu_3)$ for $\mu_1 \le \mu_2 \le \mu_3$.  We next
describe how {\boldmath$\mu$} influences the map $\phi$.

\begin{prop}
\label{strongdeg.5}
Suppose that $\phi : \p^2 \,-\!\to \p^3$ is given by $a,b,c,d$ as
above and that $\mathrm{Syz}(a,b,c,d)$ has a strong
{\boldmath$\mu$}-basis for $\mbox{\boldmath$\mu$} =
(\mu_1,\mu_2,\mu_3)$, $\mu_1+\mu_2+\mu_3 = n$.  Assume in addition
that $\phi$ is generically one-to-one and that $\mathbf{V}(a,b,c,d)
\subset \p^2$ is a local complete intersection.  Then:
\begin{enumerate}
\item The degree of the image of $\phi$ in $\p^3$ is 
\[
{\textstyle\frac12}\big(n^2 - (\mu_1^2+\mu_2^2+\mu_3^2)\big) =
\mu_1\mu_2+ \mu_1\mu_3+\mu_2\mu_3 .
\]
\item The sum of the multiplicities of the basepoints of $\phi$ is
\[
{\textstyle\frac12}\big(n^2 + (\mu_1^2+\mu_2^2+\mu_3^2)\big) = 
n^2 - (\mu_1\mu_2+ \mu_1\mu_3+\mu_2\mu_3).
\]
\end{enumerate}
\end{prop}

\prf Let $I = \langle a,b,c,d\rangle \subset R$ and $Z =
\mathbf{V}(I) \subset \p^2$.  Since $\phi$ is generically one-to-one,
we know that the image of $\phi$ is a surface in $\p^3$ of degree
\begin{equation}
\label{degreeformula}
n^2 - \sum_{p \in Z} e(\mathcal{I}_{Z,p},\mathcal{O}_{\p^2,p}),
\end{equation}
where $\mathcal{I}_{Z} \subset \mathcal{O}_{\p^2}$ is the ideal sheaf
of $Z \subset \p^2$ and $e(\mathcal{I}_{Z,p},\mathcal{O}_{\p^2,p})$ is
the multiplicity, as defined in \cite[4.5]{bh}.  (A proof of
\eqref{degreeformula} is sketched in the Appendix.)  It follows that
part 1 of the proposition is an immediate consequence of part 2.

Since $Z$ is a local complete intersection, we have
$e(\mathcal{I}_{Z,p},\mathcal{O}_{\p^2,p}) = \dim_\C
\mathcal{O}_{Z,p}$ for all $p \in Z$.  In other words,
\[
\sum_{p \in Z} e(\mathcal{I}_{Z,p},\mathcal{O}_{\p^2,p}) = \dim_\C
H^0(Z,\mathcal{O}_Z).  
\]
Using the usual vanishing theorems for the sheaf cohomology of $\p^2$,
one obtains
\[
\dim_\C H^0(Z,\mathcal{O}_Z) = \dim_\C (R/I)_d
\]
for $d \gg 0$. Since $\dim_\C R_d = \binom{d+2}{2}$, the above
resolution for $R/I$ shows that 
\begin{equation}
\label{rmodi.5}
\dim_\C (R/I)_d =
{\textstyle\frac12}\big(n^2 + (\mu_1^2+\mu_2^2+\mu_3^2)\big)
\end{equation}
for $d \gg 0$.  Combining the above equalities completes the
proof.\qed
\medskip

In using Proposition \ref{strongdeg.5}, note that 
\[
\mu_1^2+\mu_2^2+\mu_3^2 = (\mu_1-n/3)^2 + (\mu_2-n/3)^2 +
(\mu_3-n/3)^2 + n^2/3
\]
since $\mu_1+\mu_2+\mu_3 = n$.  It follows that
\[
\mu_1^2+\mu_2^2+\mu_3^2 \ge n^2/3.
\]
Let $N = \sum_{p\in Z} e(\mathcal{I}_{Z,p},\mathcal{O}_{\p^2,p})$ be
the number of basepoints, counted with multiplicity.  Combining the
above inequality with Proposition~\ref{strongdeg.5}, we see that if
$\phi$ has a strong {\boldmath$\mu$}-basis, then the number of
basepoints of $\phi$ is bounded below by
\[
N \ge 2n^2/3.
\]
Thus surface parametrizations with strong {\boldmath$\mu$}-bases have
\emph{lots} of basepoints.

The classic example of Proposition \ref{strongdeg.5} is when $n = 3$.
Then $\mu_1+\mu_2+\mu_3 = 3$ implies $\mu_1 = \mu_2 = \mu_3 = 1$
(since $\mu_i = 0$ would imply that the image of $\phi$ lies in a
plane).  It follows that the surface has degree
\[
{\textstyle\frac12}\big(3^2-(1^2+1^2+1^2)\big) = 3,
\]
and the number of basepoints is
\[
{\textstyle\frac12}\big(3^2+(1^2+1^2+1^2)\big) = 6.
\]
This, of course, is the usual representation of a cubic surface in
$\p^3$ as $\p^2$ blown up at 6 points.  In the 19th century, algebraic
geometers were aware that cubic surfaces have strong
{\boldmath$\mu$}-bases.  For example, the 1915 edition of
\cite[Vol.~II]{salmon} states on p.~264 that ``Clebsch has used the
theorem that any cubic may be generated as the locus of the
intersection of three corresponding planes, each of which passes
through a fixed point.''  This refers to three moving planes, which
Salmon called ``sheaves of planes'' (in a footnote on p.~25 of
Volume~I, Salmon notes that the term ``sheaf'' comes from the German
``B\"undel'').  We should also mention that \cite[p.~389]{sommerville}
calls a moving plane a ``bundle of planes''.  I am grateful to Tom
Sederberg for supplying these references.

One feature of Proposition \ref{strongdeg.5} is the requirement that
the basepoint locus $Z$ be a local complete intersection.  This is
because the degree of the surface naturally involves the multiplicity
of the basepoints, while the Hilbert polynomial computation given in
\eqref{rmodi.5} computes the degree of the basepoints.  As is
well-known, these agree only for a local complete intersection.  It
would be interesting to study what happens to
Proposition~\ref{strongdeg.5} when the basepoints are not a local
complete intersection.  Since having a strong {\boldmath$\mu$}-basis
is such a restrictive condition, it is possible that the basepoints
are very special.

\subsection{Syzygies Which Vanish at Basepoints} At one point in the
proof of Theorem~\ref{tensor.3}, we needed to know that a certain
syzygy on $a,b,c$ was a Koszul syzygy (this was \eqref{abcsyz.3}
$\Rightarrow$ \eqref{abckoszul.3}).  We now study what happens when
basepoints are present.

We will consider homogeneous polynomials $a,b,c \in R = \C[s,t,u]$ of
degree $n$, where $\gcd(a,b,c) = 1$.  If $a,b,c$ have no basepoints,
then they are a regular sequence, so that every syzygy
\[
Aa + Bb + Cc = 0
\]
must be a Koszul syzygy
\[
\begin{aligned}
	{A} &= \ h_1c + h_2b\\
	{B} &= -h_2a + h_3c\\
	{C} &= -h_1a - h_3b,
\end{aligned}
\]
as in \eqref{koszulsyz.3}.  

Now suppose that $a,b,c$ have some basepoints.  Then it is easy to
make examples of syzygyies which are not Koszul.  One easy
observation is that \emph{Koszul syzygies also vanish at the
basepoints}.  This leads to the following question:
\begin{equation}
\label{question.5}
\begin{aligned}
{}&\text{If $Aa+Bb+Cc = 0$ and $A,B,C$ vanish at the}\\
{}&\text{basepoints, then is $A,B,C$ a Koszul syzygy?}
\end{aligned}
\end{equation}
To make this precise, we need the following definition.

\begin{defn}
\label{vanish.5}
Let $a,b,c \in R = \C[s,t,u]$ be homogeneous with no common factors,
and let $Z = \mathbf{V}(a,b,c)$ be their basepoint locus, regarded as
a 0-dimensional subscheme of $\p^2$.  Then a homogeneous polynomial $A
\in R$ of degree $d$ {\bf vanishes at the basepoints} if $A$ is in the
kernel of the map
\[
H^0(\p^2,\mathcal{O}_{\p^2}(d)) \longrightarrow
H^0(Z,\mathcal{O}_{Z}(d)).
\]
Equivalently, $A$ vanishes at the basepoints if and only if $A$ is in
the saturation of $I = \langle a,b,c\rangle \subset R$. 
\end{defn}

It would probably be better to say ``vanishes scheme-theoretically at
the basepoints'' or ``vanishes with multiplicity at the basepoints''
in Definition~\ref{vanish.5}.  We hope that the simpler phrase
``vanishes at the basepoints'' will not cause confusion.

We say that a syzygy $Aa+Bb+Cc = 0$ \emph{vanishes at the basepoints}
if $A,B,C$ vanish at the basepoints of $a,b,c$ in the sense of
Definition~\ref{vanish.5}.  Thus \eqref{question.5} now has a precise
meaning.  

One way to think about \eqref{question.5} is that vanishing at the
basepoints is a local condition, while being a Koszul syzygy is a
global condition.  Because of this, it turns out that the answer to
\eqref{question.5} is sometimes ``no''.  For an example, suppose that
\[
\begin{aligned}
a &= s^2u + st^2\\
b &= stu + 2t^3\\
c &= t^2u+s^3.
\end{aligned}
\]
One can show without difficulty that the $p = (0,0,1) \in \p^2$ is the
unique basepoint.  It has multiplicity 4 and degree 3, and locally
looks like $\C[s,t]/\langle s^2,st,t^2\rangle$.  Using Macaulay, one
finds the syzygy of degree 5 given by
\begin{equation}
\label{counter.5}
\begin{aligned}
A &= t^2u^3-2s^2t^2u\\
B &= -stu^3+s^3tu\\
C &= st^2u^2.
\end{aligned}
\end{equation}
It is obvious that this syzygy vanishes at $p$ in the sense of
Definition~\ref{vanish.5}.  This is because every term contains either
$s^2$, $st$ or $t^2$.  With a little more work, one can show that
\eqref{counter.5} is \emph{not} a Koszul syzygy.  So
\eqref{question.5} is not always true.

So the next question is whether there exist special classes of
basepoints for which the answer to \eqref{question.5} is ``yes''.  We
do not yet have a complete answer to this question, but we do know one
class of basepoints for which this works.  

Given $a,b,c \in R$ and $Z = \mathbf{V}(a,b,c) \subset \p^2$ as usual,
we say that $p \in Z$ has \emph{embedding dimension at most one} if
the Zariski tangent space of $Z$ at $p$ has dimension $\le 1$.  One
easily sees that this is equivalent to either of the following
conditions:
\begin{itemize}
\item One of $a,b,c$ has a nonvanishing partial derivative at $p$.
\item There are local analytic coordinates $u,v$ at $p \in \p^2$ such
that near $p$, $Z$ is defined as a formal scheme by $u = v^k = 0$.  
\end{itemize}
The second characterization shows that $Z$ is a local complete
intersection at $p$.  It follows easily that $k$ is the multiplicity
of the basepoint $p$.  These basepoints were introduced under the name
\emph{aligned} by Iarrobino in 1981 \cite{iarrobino1}, though these
days, the name \emph{curvilinear} is more common---see, for example,
\cite{iarrobino2,lebarz}. 

Then we have the following result.

\begin{thm}
\label{curvilinear.5}
Suppose that $a,b,c \in R$ are as usual, and let $Z =
\mathbf{V}(a,b,c)$ be the basepoint locus.  If all basepoints are
curvilinear $($i.e., have embedding dimesion at most one\/$)$, then
every syzygy on $a,b,c$ which vanishes on $Z$ is a Koszul syzygy.
\end{thm}

Before we can give the proof, we need some preliminary results, the
first of which is the following vanishing lemma.

\begin{lem}
\label{vanishing.5}
Let $X$ be a smooth complete surface with $q := \dim_\C
H^1(X,\mathcal{O}_X) = 0$, and let $L \subset X$ be a smooth
rational curve with self-intersection $L^2 = 1$.  Then
$H^1(X,\mathcal{O}_X(mL)) = 0$ for all $m \in \Z$.
\end{lem}

\prf
Tensoring the exact sequence $0 \to \mathcal{O}_X(-L) \to
\mathcal{O}_X \to \mathcal{O}_L \to 0$ with $\mathcal{O}_X(mL)$ gives 
\[
0 \to \mathcal{O}_X((m-1)L) \to \mathcal{O}_X(mL) \to
\mathcal{O}_{\p^1}(m) \to 0
\]
since $L \simeq \p^1$ and $L^2 = 1$.  Then the long exact sequence in
cohomology yields
\begin{equation}
\label{long.5}
\begin{aligned}
{}&\qquad\qquad\quad H^0(X,\cO_X(mL)) \to H^0(\p^1,\cO_{\p^1}(m)) \to\\
{}&H^1(X,\mathcal{O}_X((m-1)L)) \to 
H^1(X,\mathcal{O}_X(mL)) \to H^1(\p^1,\mathcal{O}_{\p^1}(m)). 
\end{aligned}
\end{equation}
Using $q = 0$ and the vanishing of $H^1(\p^1,\mathcal{O}_{\p^1}(m))$
for $m \ge 0$, the lemma follows for $m \ge 0$ by induction.

Observe that $H^0(X,\cO_X(mL)) \to H^0(\p^1,\cO_{\p^1}(m))$ is an
isomorphism for $m \le 0$, since it is $\C \simeq \C$ for $m = 0$ and
$H^0(X,\cO_X(mL)) = H^0(\p^1,\cO_{\p^1}(m)) = \{0\}$ for $m < 0$.  It
follows that for $m \le 0$, \eqref{long.5} gives the exact sequence
\[
0 \to H^1(X,\mathcal{O}_X((m-1)L)) \to H^1(X,\mathcal{O}_X(mL)).
\]
Using $q = 0$ and induction, the lemma follows for $m \le 0$.\qed
\medskip

We next study curvilinear basepoints using toric geometry, as in
\cite[Chapter 2]{fulton}.  Our goal is to construct a toric blow-up of
$0 \in \C^2 = \mathrm{Spec}(\C[s,t])$ such that $I = \langle s,t^k
\rangle \subset \C[s,t]$ becomes principal on $X$.  For this purpose,
let $N = \Z^2$, with basis $e_1,e_2$.  Then let
\[
v_i = \begin{cases} ie_1+e_2 & 0 \le i \le k \\ e_1 & i = k+1
	\end{cases}
\]
and define the cones $\sigma_0,\dots,\sigma_k$ by
\[
\sigma_i = \mathrm{Cone}(v_i,v_{i+1}),\quad 0 \le i \le k.
\]
Finally, let $\Delta_k$ be the fan consisting of the $\sigma_i$
and all of their faces.  When $k = 3$, here is a picture of the fan
$\Delta_3$ in $\nr \simeq \R^2$:
\[
\begin{array}{c}
\begin{picture}(288,172)(-144,-78)
\put(-80,-60){\line(1,0){140}}
\put(-80,-60){\line(1,1){140}}
\put(-80,-60){\line(2,1){140}}
\put(-80,-60){\line(3,1){140}}
\put(-80,-60){\line(0,1){140}}
\put(-80,-20){\circle*{3}}
\put(-40,-60){\circle*{3}}
\put(-40,-20){\circle*{3}}
\put(0,-20){\circle*{3}}
\put(40,-20){\circle*{3}}
\put(-90,-15){$v_0$}
\put(-45,-15){$v_1$}
\put(-5,-15){$v_2$}
\put(35,-15){$v_3$}
\put(-45,-70){$v_{4}$}
\put(-45,40){$\sigma_0$}
\put(50,40){$\sigma_1$}
\put(50,-5){$\sigma_2$}
\put(50,-40){$\sigma_3$}
\end{picture}
\end{array}
\]

The first quadrant $\sigma =\mathrm{Cone}(e_2,e_1) =
\mathrm{Cone}(v_0,v_{k+1})$ is the union of the $\sigma_i$, and one
sees that $\Delta_k$ is obtained from $\sigma$ by a sequence of
$k$ stellar subdivisions.  Turning to the corresponding toric
varieties, $\sigma$ gives $\C^2$, and $\Delta_k$ gives the smooth
toric surface $X_k = X(\Delta_k)$.  We also have a natural map $\pi :
X_k \to \C^2$ which is the successive blow-up of smooth points.

Let $Z \subset \C^2$ be the subscheme defined by $I = \langle
s,t^k\rangle$, with ideal sheaf $\cI_Z \subset \cO_{\C^2}$.  Under $\pi :
X_k \to \C^2$, we get the \emph{inverse image ideal sheaf} 
\[
\cI'_Z = \pi^{-1}\cI_Z \cdot \cO_{X_k} \subset \cO_{X_k},
\]
as defined in \cite[II.7]{hartshorne}.  We can now explain how $\pi :
X_k \to \C^2$ relates to the ideal $I$.  Recall that $K_X$ denotes the
canonical divisor of a smooth surface $X$.

\begin{prop} 
\label{toric.5}
Let $\pi : X_k \to \C^2$ be as above.  Then $\cI'_Z = \cO_{X_k}(-E)$,
where $K_{X_k} = \pi^* K_{\C^2} + E$.
\end{prop}

\prf
As usual, each interior ray $\rho_i = \mathrm{Cone}(v_i)$, $1 \le i
\le k$, corresponds to an orbit closure $E_i \simeq \p^1$ in $X$.  We
will show that $E = E_1 + 2E_2 + 3E_3 + \dots + kE_k$ satisfies the
two conditions of the proposition.

We know that $X_k$ has the affine open cover given by
\[
U_i = \mathrm{Spec}(\C[\sigma_i^\vee\cap M]),\quad 0 \le i \le k,
\]
where $M$ is the dual of $N$.  In order to get coordinates for these
affine pieces, let $\epsilon_1,\epsilon_2$ in $M$ be the dual basis of
$e_1,e_2$ in $N$.  Then $s = \chi^{\epsilon_1}, t = \chi^{\epsilon_2}$ are
coordinates on 
\[
\C^2 = \mathrm{Spec}(\C[\sigma^\vee\cap M]) = \mathrm{Spec}(\C[s,t]),
\]
and one can show without difficulty that
\begin{equation}
\label{localcoords.5}
U_i = \begin{cases}\mathrm{Spec}(\C[t^{i+1}/s,s/t^i]) & 0
\le i \le k-1\\
	\mathrm{Spec}(\C[t,s/t^k]) & i = k.
	\end{cases}
\end{equation}

As above, the orbit closures $E_i \simeq \p^1$ correspond
to the interior rays of $\Delta_k$.  Of these, $U_0$ meets only $E_1$,
$U_1$ meets only $E_1$ and $E_2$, and so on until $U_k$ meets only
$E_k$.  Furthermore, in terms of the local coordinates
\eqref{localcoords.5}, we have equations
\[
\begin{aligned}
\text{On}\ U_0, &\ E_1\ \text{is defined by}\ s = 0\\
\text{On}\ U_1, &\begin{cases} E_1\ \text{is defined by}\ t^2/s = 0 &
		\\[-2pt] 
		E_2\ \text{is defined by}\ s/t = 0 & \end{cases}\\[-2pt] 
&\ \, \vdots \\[-1pt] 
\text{On}\ U_{k-1}, &\begin{cases} E_{k-1}\ \text{is defined by}\
		t^k/s = 0 & \\[-2pt] 
		E_k\ \text{is defined by}\ s/t^{k-1} = 0 &
		\end{cases}\\[3pt] 
\text{On}\ U_k, &\ E_k\ \text{is defined by}\ s = 0.
\end{aligned}
\]

Let $E = E_1 + 2E_2 + 3E_3 + \dots + kE_k$ be as above.  Using the
local equations for the $E_i$, one can verify that
\begin{equation}
\label{eeq.5}
\text{On}\ U_i, \ E\ \text{is defined by} \begin{cases} s= 0 & 0 \le
i \le k-1 \\ t^k = 0 & i = k.\end{cases}
\end{equation}
It remains to show that this divisor has the desired properties.

First, the inverse image ideal sheaf $\cI'_Z$, when restricted to
$U_i$, corresponds to the ideal generated by $s,t^k$ in the coordinate
ring of $U_i$.  Using \eqref{localcoords.5}, one easily computes that
in the coordinate ring of $U_i$, one has
\[
\langle s,t^k\rangle = \begin{cases} \langle s\rangle & 0 \le
i \le k-1 \\ \langle t^k\rangle & i = k.\end{cases}
\]
By \eqref{eeq.5}, this is the ideal defining $E$ on $U_i$, which
implies that
\[
\cI'_Z = \cO_{X_k}(-E).
\]

Second, we need to show that $K_{X_k} = \pi^* K_{\C^2} + E$.  This
follows by representing $X_k \to \C^2$ as a composition of successive
blow-ups of smooth points 
\[
X_k \xrightarrow{\ \pi_k\ } X_{k-1} \xrightarrow{\ \pi_{k-1}\ } \cdots
\xrightarrow{\ \pi_2\ } X_1 \xrightarrow{\ \pi_1\ } X_0 = \C^2
\]
and using the well-known fact that $K_{X_{i}} = \pi_i^* K_{X_{i-1}} +
E_i$.  The key point that needs to be checked is that $\pi_{i+1}$
blows up a point on $E_{i}\setminus (E_1\cup\dots\cup E_{i-1})$.  We
omit the straightforward details.\qed
\medskip

We can now prove Theorem \ref{curvilinear.5}.
\medskip

\noindent {\bf Proof of Theorem \ref{curvilinear.5}.\ \ }
Since the basepoint locus $Z = \mathbf{V}(a,b,c)$ is curvilinear, each
$p \in Z \subset \p^2$ is analytically equivalent to the germ of the
analytic space $0 \in \mathbf{V}(s,t^k) \subset \C^2$, where $k$
depends on $p$.  By Proposition~\ref{toric.5}, there is a smooth
surface $\pi : X \to \p^2$ such that $\cI_Z' = \cO_X(-E)$ and $K_X =
\pi^* K_{\p^2} + E$.

By the definition of $Z$, $a,b,c$ give global sections of $\cI_Z(n)$.
In the blow-up, these become sections $\tilde a,\tilde b,\tilde c$
(the \emph{proper transforms} of $a,b,c$) of $\cI_Z'\otimes_{\cO_X}\!
\pi^*\cO_{\p^2}(n)$.  If we let $L\subset X$ be the inverse image of a
line in $\p^2$ missing $Z$, then this implies that $\tilde a,\tilde
b,\tilde c$ are global sections of $\cO_X(nL-E)$.  The key point is
that $\tilde a,\tilde b,\tilde c$ have no basepoints on $X$.  This, of
course, is why we blew-up $\p^2$.  Recall why this is true:\ $a,b,c$
give an exact sequence
\[
\cO_{\p^2}(-n)^3 \to \cI_Z \to 0,
\]
so that on the blow-up, 
\[
\cO_X(-nL)^3 \to \cI_Z' \to 0
\]
is still exact.  Since $\cI_Z' = \cO_X(-E)$ is invertible, $\tilde
a,\tilde b,\tilde c$ give the exact sequence
\begin{equation}
\label{beginkos.5}
\cO_X(E-nL)^3 \to \cO_X \to 0.
\end{equation}
Hence these sections can't vanish simultaneously on $X$.  Once we know
that $\tilde a,\tilde b,\tilde c$ have no basepoints, they are locally
a regular sequence, so that we can extend \eqref{beginkos.5} to the
Koszul complex
\begin{equation}
\label{koszul1.5}
0 \to \cO_X(3E-3nL) \to \cO_X(2E-2nL)^3 \to \cO_X(E-nL)^3 \to \cO_X
\to 0,
\end{equation}
which is exact on $X$.

Now suppose that $Aa + Bb + Cc = 0$ is a syzygy, where $A,B,C$ have degree
$m$ and vanish at $Z$.  This means that $A,B,C$ are global sections of
$\cI_Z(m)$, and taking proper transforms as above, we get global
sections $\widetilde A,\widetilde B,\widetilde C$ of $\cO_X(mL-E)$.
Now tensor \eqref{koszul1.5} with $\cO_X((m+n)L-2E)$ to obtain
\[
\begin{aligned}
0 \to \cO_X((m-2n)L+E) &\to \cO_X((m-n)L)^3 \to \\
&\cO_X(mL-E)^3 \to \cO_X((m+n)L-2E) \to 0.
\end{aligned}
\]
Splitting this into two short exact sequences in the usual way, one
sees that if
\begin{equation}
\label{vanishing2.5}
H^1(X,\cO_X((m-2n)L+E)) = \{0\},
\end{equation}
then
\begin{equation}
\label{exactness.5}
\begin{aligned}
{}&H^0(X,\cO_X((m-n)L))^3 \to \\
{}&\qquad\qquad\qquad H^0(X,\cO_X(mL-E))^3 \to H^0(X,\cO_X((m+n)L-2E)) 
\end{aligned}
\end{equation}
is exact at the middle term.  For more details, see Proposition~2.2 of
\cite{cgz}.

We claim that \eqref{vanishing2.5} holds for all $m$.  To see why,
note that $K_X = \pi^* K_{\p^2} + E = -3L+E$.  Thus Serre duality
implies that
\[
\begin{aligned}
H^1(X,\cO_X((m-2n)L+E)) &\simeq H^1(X,\cO_X(-[(m-2n)L+E] + K_X))^*\\
	 &\simeq H^1(X,\cO_X((2n-m-3)L))^*.
\end{aligned}
\]
However, Lemma \ref{vanishing.5} implies $H^1(X,\cO_X((2n-m-3)L)) =
\{0\}$, and \eqref{vanishing2.5} follows.  

This shows that \eqref{exactness.5} is exact at the middle term.
However, the syzygy $Aa+Bb+Cc = 0$ implies that $(\widetilde
A,\widetilde B,\widetilde C) \in H^0(X,\cO_X(mL-E))^3$ maps to 0 in
$H^0(X,\cO_X((m+n)L-2E))$.  By exactness, this means that there is
$(h_1,h_2,h_3) \in H^0(X,\cO_X((m-n)L))^3$ which maps to $(\widetilde
A,\widetilde B,\widetilde C)$ in \eqref{exactness.5}.  However,
\[
H^0(X,\cO_X((m-n)L)) \simeq H^0(\p^2,\cO_{\p^2}(m-n))
\]
since $\pi_*\cO_X = \cO_{\p^2}$.  It follows that $(h_1,h_2,h_3)$
are homogeneous polynomials of degree $m-n$ which make $A,B,C$ into a
Koszul syzygy, as claimed.\qed
\medskip

The proof of this theorem shows that a curvilinear basepoint locus $Z
\subset \p^2$ has the property that there is a blow-up $\pi : X \to
\p^2$ such that $\cI_Z' = \cO_X(-E)$ and $K_X = \pi^* K_{\p^2} +E$.
Obviously, we can replace $\p^2$ with any smooth complete surface and
the same result holds.  Furthermore, we suspect that the converse is
true, though we have not proved this.

The main unresolved question is whether curvilinear basepoints are the
only basepoints which have the property of
Theorem~\ref{curvilinear.5}.  For example, does the theorem hold when
$Z$ is a local complete intersection?  When blowing up more
complicated basepoints, Enriques introduced \emph{Enriques diagrams}
to keep track of the combinatorics of the blow-ups.  This has been
greatly generalized in \cite{cgs}.  Our treatment of curvilinear
basepoints uses toric geometry and is a special case of the
\emph{toric clusters} and \emph{toric constellations} discussed in
\cite{gsp}.

It would also be interesting to explore higher-dimensional versions of
Theorem~\ref{curvilinear.5}.  It is possible that the papers just
mentioned might provide useful tools for attacking this problem.

\section{Final Comments} 

The results in the paper raise many
questions to pursue.  For instance, the results of Sections 3 and 4
should extend to the case when there are basepoints.  To give an
example of how this might work, suppose that
\[
\phi : \p^1\times\p^1 \,-\!\to \p^3
\]
has bidegree $(m,n)$ with a single base point $p$ of multiplicity one.
Then the map $MP$ of \eqref{mp.3} cannot be onto, since its image lies
in the subspace of polynomials of degree $2n-1$ which vanish at $p$.
Here is what we would expect to happen:
\begin{itemize}

\item The basepoint should cause the rank of $MP$ should drop by one,
so that there should be exactly one linearly independent moving plane
of degree $(m-1,n-1)$ which follows $\phi$.

\smallskip

\item The basepoint should also cause the rank of $MQ$ from
\eqref{mq.3} to drop by three.  This makes sense since
$a^2,ab,\dots,cd,d^2$ all vanish to second order at $p$.  This means
there should be $mn+3$ linearly independent moving quadrics of degree
$(m-1,n-1)$ which follow $\phi$.

\smallskip

\item Multiplying the moving plane of the first bullet by $x,y,z,w$ as
we did in Section 4 gives four moving quadrics.  Picking a
complementary subspace, we get $mn-1$ linearly independent moving
quadrics which don't come from moving planes.  

\smallskip

\item If we use the one moving plane and $mn-1$ moving quadrics to
create a matrix $M$, we get an $mn\times mn$ matrix with one linear
row and $mn-1$ quadratic rows.  It follows that $\det M$ has degree $1
+ 2\cdot(mn-1) = 2mn-1$.

\smallskip

\item Since the image of $\phi$ has degree $2mn-1$ (assuming $\phi$ is
generically one-to-one), the equation of the image should be $\det M =
0$.
\end{itemize}
As the number and complexity of the basepoints increases, it becomes
less clear how to modify the matrix $M$ in order to get the equation
of the image.  

At the end of Section 5, I mentioned some open questions concerning
syzygies and basepoints.  There are also numerous questions from
Sections 3 and 4 that I would love to be able to answer, including the
following:
\begin{itemize}

\item For the tensor product case, we used syzygies of bidegree
$(m-1,n-1)$, and for the triangular case, we used degree $n-1$.
What is the systematic reason for choosing these degrees?
\item Why do we need $a^2,ab,\dots,cd,d^2$ in order to compute the
implicit equation? In a sense, we can think of $a^2,ab,\dots,cd,d^2$
as the ``second symmetric power'' of $a,b,c,d$.  

\smallskip

\item When $I = \langle a,b,c,d\rangle \subset R = \C[s,t,u]$ is the
ideal coming from a triangular surface parametrization, what is the
free resolution of $R/I$?  How do the Betti numbers depend on $n$?  As
far as I know, this is an open problem.

\smallskip

\item Resolutions of length 3 are studied in \cite{weyman}.  Do the
results of this paper shed any light on the resolution of $R/I$?  One
intriguing observation is that symmetric powers (such as those
mentioned in the second bullet) appear in Weyman's description.

\smallskip

\item In the curve case, the determinant giving the implicit equation
is the resultant of a basis of the syzygy module.  In the surface
case, $\mathrm{Syz}(a,b,c,d)$ is usually not free.  But is it still
possible to interpret the determinants of Theorems~\ref{tensor.3}
and~\ref{triangular.4} as some sort of resultant?  Klaus Altmann
suggested that this may involve the determinant of a complex built
from a free resolution of $\mathrm{Syz}(a,b,c,d)$ (which is related to
the free resolution of $R/I$ mentioned in the previous two bullets).  

\smallskip

\item In the curve case, we saw that the structure of the syzygy
module was closely related to the regularity of $R/I$.  Does something
similar happen in the surface case?  Also, how does one define the
regularity of a tensor product surface?  This is not obvious since
everything is now bigraded.

\smallskip

\item Beauville studies hypersurfaces in $\p^n$ whose defining
equation is a determinant in \cite{beauville}.   How do his results
relate to the theorems of Sections~3 and~4? 

\smallskip

\item When the base points are a complete intersection, the recent
preprint \cite{buse} shows how the implicitization problem can be
solved using the \emph{residual resultants} defined in \cite{bem}.  In
\cite{buse}, these resultants are computed as the gcd of three
determinants or the product of two determinants divided by a third.
One question is whether one can combine residual resultants with the
syzygy methods introduced here to give a \emph{determinantal} formula
for the implicit equation, assuming that the base points are a local
complete intersection.  For example, in the tensor product case, we
explained at the beginning of this section how to modify the matrix
$M$ in the presence of a single basepoint.  Can the determinant of $M$
be interpreted as a residual resultant?
\end{itemize}

\section{Acknowledgements} 

First, I would like to thank my coauthors Tom Sederberg, Falai Chen,
Ron Goldman and Ming Zhang for asking interesting questions in
algebraic geometry and commutative algebra.  I am also grateful to
Carlos D'Andrea, Klaus Altmann, Laurent Bus\'e, Hal Schenk and Jerzy
Weyman for useful comments, and I would like to thank Tony Iarrobino
for telling me about aligned and curvilinear singularities.  Thanks
also go to the referees for many helpful suggestions.

A preliminary version of this article was presented at the
AMS-IMS-SIAM summer conference at Mount Holyoke organized by Ed Green,
Serkan Hosten, Reinhard Laubenbacher and Vicki Powers.  Also, portions
of Sections 3, 4 and 5 were presented at a Special Session organized
by Irena Peeva at the AMS meeting in Lowell.  Thanks to all of these
organizers for their hard work.

\section*{Appendix}

In the proof of Proposition~\ref{strongdeg.5}, we gave a formula
\eqref{degreeformula} for the degree of the image of a generically
one-to-one map $\phi : \p^2 \,-\!\to \p^3$.  As I learned by reading
\cite{buse}, this formula follows easily from the results in
\cite{fulton1}.

Let $\phi : \p^2 \,-\!\to \p^3$ be defined by homogeneous
polynomials $a,b,c,d \in \C[s,t,u]$ of degree $n$.  As in Section 5,
the set of basepoints is $Z = \mathbf{V}(a,b,c,d) \subset \p^2$.  We
will assume that:
\begin{itemize}
\item $\gcd(a,b,c,d) = 1$, so that $Z$ is finite.  The ideal sheaf of
$Z$ is $\mathcal{I}_Z \subset \mathcal{O}_{\p^2}$ and, for each $p \in
Z$, we get the multiplicity $e(\mathcal{I}_{Z,p},\mathcal{O}_{\p^2,p})$.
\item $S = \overline{\phi(\p^2 \setminus Z)}$ is a surface in $\p^3$.
This allows us to define the degree $\deg(S)$ of $S \subset \p^3$ and
the generic degree $\deg(\phi)$ of $\phi$.  
\end{itemize}
These numbers are related by the \emph{degree formula}, which goes as
follows:
\[
\deg(\phi)\,\deg(S) = n^2 -  \sum_{p \in Z}
e(\mathcal{I}_{Z,p},\mathcal{O}_{\p^2,p}).
\]
To prove this, set $L = \mathcal{O}_{\p^2}(n)$ and let $\pi : X \to
\p^2$ be the blow-up of $\mathcal{I}_{Z} \subset \mathcal{O}_{\p^2}$.
Since the cycle $\pi_*[X]$ has degree $\deg(\phi)\,\deg(S)$,
\cite[Proposition~4.4]{fulton1} implies that
\[
\deg(\phi)\,\deg(S) = \int_{\p^2} c_1(L)^2 - \int_Z (1+c_1(L))^2 \cap
s(Z,\p^2),
\]
where $s(Z,\p^2)$ is the Segre class.  It is standard that
$\int_{\p^2} c_1(L)^2 = n^2$.  Furthermore, $c_1(L) = 0$ on $Z$ 
since $Z$ is zero-dimensional.  Hence the above formula 
reduces to
\[ 
\deg(\phi)\,\deg(S) = n^2 - \int_Z 1 \cap s(Z,\p^2).
\]
However, as explained on \cite[p.~79]{fulton1}, the Segre class
$s(Z,\p^2)$ is given by
\[
s(Z,\p^2) = \sum_{p\in Z} e(\mathcal{I}_{Z,p},\mathcal{O}_{\p^2,p})
[p].
\]
The degree formula now follows immediately.  

When $\phi$ is generically one-to-one, the degree formula reduces to
the formula \eqref{degreeformula} used in the proof of
Theorem~\ref{strongdeg.5}.

\end{document}